\documentclass[12pt]{article}
\usepackage{lineno}
\modulolinenumbers[5]
\usepackage{caption}
\usepackage[a4paper]{geometry}
\usepackage{amsfonts}
\usepackage{arydshln}
\usepackage{xtab}
\usepackage{setspace}
\usepackage{fancyhdr}
\usepackage[dvips]{graphicx}
\usepackage[T1]{fontenc}
\usepackage[latin1,utf8]{inputenc}
\usepackage[main=english,italian]{babel}
\usepackage{type1ec}
\usepackage[final]{microtype}
\usepackage{amsmath,amssymb,amsthm,eucal,dsfont,bm,mathrsfs,stmaryrd}
\usepackage{lmodern}
\usepackage{textcomp}
\usepackage{pict2e,enumitem}
\usepackage{hyperref}
\usepackage[dvipsnames]{xcolor}
\usepackage[nodayofweek]{datetime} % change the format of printed dates (no american style!)
\usepackage{comment}

\hypersetup{
    pdfpagemode={UseOutlines},
    bookmarksopen,
    pdfstartview={FitH},
    colorlinks,
    linkcolor={black},
    citecolor={black},
    urlcolor={black}
            }

\definecolor{cornellred}{rgb}{0.7, 0.11, 0.11}

\usepackage{geometry}

\newdateformat{monthyear}{\monthname[\THEMONTH], \THEYEAR}% new date format

\theoremstyle{plain}

\theoremstyle{definition}

\renewcommand
        {\thefootnote}{\arabic{footnote}}
        
\newcommand{\symfootnote}[1]{%
\let\oldthefootnote=\thefootnote%
\stepcounter{mpfootnote}%
\addtocounter{footnote}{-1}%
\renewcommand{\thefootnote}{\fnsymbol{mpfootnote}}%
\footnote{#1}%
\let\thefootnote=\oldthefootnote%
}

\usepackage{tikz}

\newcommand\blfootnote[1]{%
  \begingroup
  \renewcommand\thefootnote{}\footnote{#1}%
  \addtocounter{footnote}{-1}%
  \endgroup
}

\makeatletter
\def\bbibitem#1{\item[]%
    \if@filesw\immediate\write\@auxout{\string \bibcite {#1}{\the\value{\@listctr }}}\fi\ignorespaces}
\makeatletter
\title{Reasoning by Analogy in Mathematical Practice}
\author{Nicolò Cangiotti\footnote{Politecnico di Milano, Department of Mathematics, via Bonardi 9, Campus Leonardo, 20133, Milan (Italy). E-mail: \texttt{nicolo.cangiotti@polimi.it}} \ \& Francesco Nappo\footnote{Politecnico di Milano, Department of Mathematics, via Bonardi 9, Campus Leonardo, 20133, Milan (Italy). E-mail: \texttt{francesco.nappo@polimi.it}}}
\date{}

\begin{document}

\maketitle

\begin{abstract}
The testimony and practice of notable mathematicians indicate that there is an important phenomenological and epistemological difference between superficial and deep analogies in mathematics. In this paper, we offer a descriptive theory of analogical reasoning in mathematics, stating general conditions under which an analogy may provide genuine inductive support to a mathematical conjecture (over and above fulfilling the merely heuristic role of ‘suggesting’ a conjecture in the psychological sense). The proposed conditions generalize the criteria put forward by Hesse (1963) in her influential work on analogical reasoning in the empirical sciences. By reference to several case-studies, we argue that the account proposed in this paper does a better job in vindicating the use of analogical inference in mathematics than the prominent alternative defended by Bartha (2009). Moreover, our proposal offers novel insights into the practice of extending to the infinite case mathematical properties known to hold in finite domains.
\blfootnote{Both authors contributed equally to this research.}
\end{abstract}

\section{Introduction}
\label{Sec1}
How do mathematicians make mathematical discoveries? Any historically accurate answer to this question cannot ignore the role of analogies. Explaining his efforts in the theory of fields, notable mathematician André Weil writes to his sister Simone (the famous philosopher): “Artin functions \emph{in the abelian case} are \emph{polynomials}, which one can express by saying that these fields furnish a \emph{simplified} model of what happens in number fields; here, there is thus room to conjecture that the non-abelian Artin functions are still polynomials; \emph{that is just what occupies me at the moment}” (in Krieger 2003:340). As Weil elaborates, this way of making conjectures about the features of an unfamiliar domain from those of a more well-understood analogue is a rather indispensable part of a mathematician’s methodological toolbox. Although the “faculty of seeing connections among things that are in appearance completely different... does not fail to lead us astray quite often” (340), little progress in the history of mathematics would have been achieved without an attentive study of the relations that distinct mathematical domains entertain with one another - “their conflicts and delicious reciprocal reflections, their furtive caresses, their inexplicable quarrels” (339). As a matter of historical fact, “nothing is more fecund than these slightly adulterous relationships; and nothing gives greater pleasure to the connoisseur” (339).

From a methodological perspective, it seems natural to equate the role of analogy in mathematical research to that of a useful heuristics - a psychological aid for generating conjectures to be proved or disproved. However, a closer look reveals that the attribution of a heuristic role insufficiently captures what analogies mean for the working mathematician. Identifying an analogy with a well-understood mathematical domain is often the basis of considerable expectations about the truth of a conjecture. When the aim is to devise a proof, entire research programs are initiated and pursued that exploit analogies with mathematical problems already solved. These aspects of mathematical practice indicate an implicit commitment to the view that drawing an analogy with a familiar mathematical domain can not only be useful heuristically, but also be a source of \emph{inductive support} for yet unproven conjectures, i.e., some (possibly small, but non-negligible) additional ground for believing them.\footnote{Poincaré (1905) is an eminent early source for a defense of the inductive role of analogy in mathematical research. We take no stance in what follows as to how ‘inductive support’ should be cashed out - if, for instance, it should be understood in terms of \emph{acceptability} (Bartha 2009), or as a probabilistic notion, or a mixture; see also section seven.}

We find a contemporary illustration of the inductive role of analogy in the history of the Poincaré conjecture. Formulated by Poincaré  in 1904, this theorem of algebraic geometry states that every simply connected, closed $3$-manifold (a generalization of Euclidean space) is homeomorphic to a $3$-sphere (a higher-dimensional analogue of a sphere). While resisting several attempts at a proof, the conjecture was made plausible by an analogue theorem in the two-dimensional case, viz., that every simply connected, compact $2$-dimensional surface without boundaries is homeomorphic to the $2$-sphere (the ‘ordinary’ sphere). The proof-idea was later found with the help of another analogy, with the so-called ‘Ricci flow with surgery’. In brief, Hamilton (1982) noted that one could perform operations on $3$-manifolds that are mathematically analogous to the action of the heat equation; he thus proposed that the effect of a Ricci flow on the former tends to a behavior analogous to the uniform behavior that the heat equation tends to under the same operation. Eventually, Perel’man (2003) closed the remaining gap in Hamilton’s proof-idea, at a time when practically no more doubt remained about the conjecture’s truth.

Given the centrality of analogy to mathematics, it is surprising that so little work in contemporary philosophy has been devoted to the topic.\footnote{Relevant works include Grosholtz (1985), Pickering (1997), Corfield (2008), Bartha (2009), and Mancosu (2009).} In particular, marginal attention has been given to the problem of articulating credible foundations for the inductive use of analogy in mathematics. Yet, the challenge is arguably no less pressing for the realm of mathematics as it is for the domain of the empirical sciences. Perhaps some insight into the problem may be gained by considering what Hacking (2014) has called the “astonishing” (6) fact about mathematics- the fundamental analogy between geometry, algebra and arithmetic. This extraordinary concurrence of branches otherwise so diverse from one another, both in subject matter and historical origin, is echoed in André Weil’s remark that: “my work consists in deciphering a trilingual text (cf. the Rosetta stone); of each of the three columns I have only disparate fragments; I have some ideas about each of the three languages; but I know as well that there are great differences in meaning from one column to another” (340). Exactly how to turn these vague hints (however suggestive) into a solution to the problem of justification remains a difficult and open question. 

Our aim in this paper is to make progress on a related philosophical issue, which can be regarded as propaedeutic to addressing the problem of justification. This is the broadly ‘descriptive’ project to provide an account of the conditions under which an analogy in mathematics is considered capable of non-negligible support to a mathematical claim - quite regardless of the ultimate stance one assumes on the problem of justification.\footnote{We regard this as a contribution to the “philosophy of mathematical practice” (Mancosu 2008). For introductions to the topic of non-deductive reasoning in mathematical research, see Polya (1954a, 1954b) and  Franklin (2013).} The conditions that we will put forward in this paper are meant to apply both to analogical arguments used in support of a mathematical conjecture and those used in support of some undemonstrated claim which helps prove a conjecture. We will defend our framework for distinguishing strong from weak analogical inferences in mathematics by reference to accessible historical case-studies. One of the main promises of this account will be its potential to clarify how analogical reasoning can be used to extend mathematical properties known to hold in the finite case to infinite domains. 

An important theme in the following sections will be elaborating the divergences between the account proposed in this paper and Bartha’s (2009) - arguably the most developed attempt at an account of plausible analogical inference in the current philosophical literature. As we will argue, Bartha is too quick to abandon the ‘two-dimensional’ framework for analogical reasoning in empirical sciences put forward by Hesse (1963). Contra Bartha’s claim that “in mathematics... material analogy, in Hesse’s sense, plays no role at all” (43), we will defend the necessity of restricting inductively significant analogical inferences in mathematics to those which rely upon ‘\emph{material}’ or ‘\emph{pre-theoretic}’ similarities in Hesse’s sense. Furthermore, we will argue that, corresponding to Hesse’s distinction between analogical arguments in science that project \emph{causal relations} from source to target and those that project merely \emph{statistical correlations}, there exists a parallel need to restrict inductively significant analogical inferences in mathematics to only those that project some \emph{robust mathematical connection} from source to target domain. The proposed case-studies will aim to expose these dramatic divergences in a clear and accessible manner. 

\medskip

The discussion below will proceed as follows. In section two, we will consider some illustrative case-studies of analogical reasoning at work in mathematical research. The discussion will help clarify how mathematical investigation by analogy looks from the perspective of a working mathematician. Borrowing two expressions from Hadamard (1954), we will distinguish mere analogical “hookings”, which play mostly a heuristic role, and analogical “relay-results”, which instead are signs of a deeper and more telling similarity. In section three, we will contrast the prominent descriptive approaches by Hesse (1963) and Bartha (2009). In sections four to six, we will articulate a novel account that, generalizing Hesse’s (1963) framework,  aims to capture what makes a given similarity (or set thereof) an analogical ‘relay-result’ in our sense. We will show the strengths of the proposed conditions vis-a-vis Bartha’s (2009) by illustrating their application to several case-studies, taken from diverse areas of mathematics. Finally, section seven will conclude with a summary of our novel account and with directions for future research.

\section{Phenomenology of Analogical Discovery}
\label{Sec2}
In this section, our aim is to make explicit some of the intuitive judgments that accompany the working practice of mathematicians when reasoning by analogy. To avoid making the discussion unnecessarily complex, we will provide adaptations of actual historical examples to bring out features of the mathematical mind in a way that can be accessible to non-mathematicians. The idea that we will aim to illustrate is that a central part of mathematical training consists in developing a quasi-perceptual capacity to distinguish superficial from deeper analogies. We think of this capacity as an instance of Hadamard’s (1954) discriminating faculty of the mathematical mind, whereby one is led to recognize a difference between mere ‘hookings’ (77) and ‘relay-results’ (80). The discussion below will be useful for motivating the account to be presented in section three, where the phenomenologically salient notions of depth and superficiality will be explicated in a rigorous way and related to the notion of strength of an analogical argument. 
\subsection*{\normalsize{2.1 Euler Characteristic}}
Let’s start with a classic case of analogy in mathematics (cf. Lakatos 1976 and Polya 1954a for further takes on this suggestive historical example). It starts from the observation that, in plane geometry, any convex polygon has exactly the same number of vertices ($V$) and edges ($E$):  
\begin{equation}
\label{e1}
    V=E.
\end{equation}
The question Euler (1758) posed was whether any analogous regularity holds for convex polyhedra. His conjecture, which was eventually proved by Cauchy (1813), is known as the \emph{Euler characteristic}. It states that, if $F$ is the number of faces of a polyhedron, then:
\begin{equation}
\label{e2}
V+F=E+2.
\end{equation}
For the sake of the following discussion, let’s pretend that we did not know the answer to Euler’s question. How would a trained mathematician reason to the formula for the three-dimensional case? For a start, equation \eqref{e1} could induce a direct form of analogy, as in: 
\begin{equation}
\label{e3}
    E=F.
\end{equation}
The conjecture would exploit the fact that, since an extra dimension is added when moving from plane to solid geometry, edges (lines) in solids may correspond to vertices (points) in a plane and faces (surfaces) in solids may correspond to edges (lines) in a plane. This mode of reasoning by proportion is common in geometry (see also Example 4.2) and could make it plausible that convex polyhedra obey equation \eqref{e3}. As further reflection shows, however, equation \eqref{e3} is false: a cube has twelve edges but only six faces. In this instance, the analogy leads us astray.

\begin{figure}[ht!]
\centering
\begin{tikzpicture}[scale=1,>=latex]
 %\draw[help lines, step=0.5] (-3,0) grid (8,5);
  \fill (0,0) circle[radius=3pt];
  \fill (5,0) circle[radius=3pt];
  \fill (5,5) circle[radius=3pt];
  \fill (0,5) circle[radius=3pt];
  \fill (1.5,1.5) circle[radius=3pt];
  \fill (3.5,1.5) circle[radius=3pt];
  \fill (3.5,3.5) circle[radius=3pt];
  \fill (1.5,3.5) circle[radius=3pt];
%%%%%%%%%%%%%%%%%%%%
  \draw[ultra thick] (0,0)  rectangle (5,5);
  \draw[ultra thick]  (1.5,1.5) rectangle (3.5,3.5);
  \draw[ultra thick] (0,0)--(1.5,1.5);
  \draw[ultra thick] (5,5)--(3.5,3.5);
  \draw[ultra thick] (5,0)--(3.5,1.5);
  \draw[ultra thick] (0,5)--(1.5,3.5);
  %\draw[ultra thick, densely dotted] (0,0) -- (0,6);
%NODI
\end{tikzpicture}
{\caption*{Figure 1}}
\end{figure}

\bigskip

Let’s consider a different attempt, which is adapted from Polya’s (1954:43) discussion. After some tedious counting and checking, one notes that cubes, tetrahedra and dodecahedra all satisfy the equation (where ‘$S$’ stands for the number of ‘solids’, or three-dimensional elements): 
\begin{equation}
\label{e4}
V-E+F-S=1.
\end{equation}
For instance, a cube satisfies \eqref{e4} since it has eight vertices, twelve edges, six faces, and one solid; hence, the alternating sum is one. Additionally, we may note that \eqref{e2} can be rewritten as: 
\begin{equation}
\label{e5}
V-E+F=1 .
\end{equation}
While we have only checked \eqref{e4} in few instances, the similarities with \eqref{e5} suggest that the regularity we are looking after is expressed by an alternating sum of the number of zero-dimensional entities (vertices), one-dimensional entities (edges), two-dimensional entities (faces), three-dimensional entities (solids), and so on. If this is correct, \eqref{e4} is a true conjecture.

It turns out that the conclusion that this second reasoning leads to is correct: trivially, in three dimensions \eqref{e4} is equivalent to \eqref{e2}. But assuming that we do not yet have a proof of Euler’s characteristic, one might ask: how strong is the argument from the algebraic resemblance between \eqref{e4} and \eqref{e5}? In particular, does the resemblance provide non-negligible additional ground for believing \eqref{e4}? To address this question, let’s compare the algebraic argument above to a third way of reasoning about Euler’s problem, which again involves the use of analogy.

Imagine taking any polyhedron and compressing it until it is flat. The edges and faces of the polyhedron cannot break, but they can (and most likely will) get stretched and deformed. The resulting, smashed polyhedron has its vertices connecting exactly the same edges and faces as before, but is now a two-dimensional figure. As Fig. 1 illustrates for the cube, this operation leaves both $E$ and $V$ unchanged, but in the compression we lose one face ($F$). Counting the smashed cube’s elements, and adding the missing face, we note that \eqref{e4} holds for the cube.

\begin{figure}[ht]
\centering
\begin{minipage}[c]{.49\textwidth}
\centering
\begin{tikzpicture}[scale=1,>=latex]
 %\draw[help lines, step=0.5] (-3,0) grid (8,5);
%%%%%%%%%%%%%%%%%%%%%%
  \draw[ultra thick] (0,0)  -- (5,0);
  \draw[ultra thick] (0,0)  -- (0,5);
  \fill[fill=black,opacity=0.75] (1.5,3.5)--(5,5)--(3.5,3.5);
  \draw[ultra thick]  (1.5,1.5) rectangle (3.5,3.5);
  \draw[ultra thick] (0,0)--(1.5,1.5);
  %\draw[ultra thick] (5,5)--(3.5,3.5);
  \draw[ultra thick] (5,0)--(3.5,1.5);
  \draw[ultra thick] (0,5)--(1.5,3.5);
  \draw[ultra thick, dotted] (1.5,1.5) -- (3.5,3.5);
  \draw[ultra thick, dotted] (1.5,1.5) -- (0,5);
  \draw[ultra thick, dotted] (0,0) -- (3.5,1.5);
\end{tikzpicture}
  \caption*{Figure 2}
\end{minipage}
\begin{minipage}[c]{.49\textwidth}
\centering
\begin{tikzpicture}[scale=1,>=latex]
 %\draw[help lines, step=0.5] (-3,0) grid (8,5);
%%%%%%%%%%%%%%%%%%%%%%
  \draw[ultra thick] (0,0)  rectangle (5,5);
  \fill[fill=black,opacity=0.75] (3.5,3.5)--(5,5)--(5,0);
  \draw[ultra thick]  (1.5,1.5) rectangle (3.5,3.5);
  \draw[ultra thick] (0,0)--(1.5,1.5);
  %\draw[ultra thick] (5,5)--(3.5,3.5);
  \draw[ultra thick] (5,0)--(3.5,1.5);
  \draw[ultra thick] (0,5)--(1.5,3.5);
  \draw[ultra thick, dotted] (1.5,1.5) -- (3.5,3.5);
  \draw[ultra thick, dotted] (1.5,1.5) -- (0,5);
  \draw[ultra thick, dotted] (0,0) -- (3.5,1.5);
  \draw[ultra thick, dotted] (1.5,3.5) -- (5,5);
%NODI
\end{tikzpicture}
\caption*{Figure 3}
\end{minipage}
\end{figure}

The next step is crucial. We imagine drawing (roughly) diagonal lines among the unconnected vertices of the smashed polyhedron to divide the latter up in (roughly) triangles.  In drawing these diagonals, both $V$ and $E$ increase by one, so $V-E$ remains constant. If we then consider a smashed cube, tetrahedron and dodecahedron, we note that the removal of any (roughly) triangular surface involves one of the following two options: either an edge and a face is removed (Fig. 2), or a vertex, two edges and a face (Fig. 3). In either case, $V-E+F$ remains constant. It follows that, if the smashed polyhedron under consideration satisfied \eqref{e4}, so would the same smashed polyhedron without one of the (roughly) triangular surfaces. Moreover, if we were to remove all (rough) triangles except one, the resulting figure would still satisfy \eqref{e4}.

Having verified that the operations of compression and triangularization yield \eqref{e4} as the correct regularity for cubes, tetrahedra and dodecahedra, we now reason - by analogy - that the same operations will also yield \eqref{e4} when considering a generic convex polyhedron. The reasoning is analogical in that it derives its plausibility, not so much from the \emph{number} of observed instances (which are relatively few), but from the \emph{similarities} that link all convex polyhedra: specifically, the fact that, because of the underlying commonalities that make them instances of the same geometrical genus, they all appear to be amenable to the same compressing and triangulating treatment to which we have subjected our cubes, tetrahedra and dodecahedra.

The new reasoning supports exactly the same (correct) conclusion as the algebraic, Polya-inspired argument considered earlier. However, it seems that the former is a much stronger reasoning than the latter\footnote{\label{fnLak}Lakatos (1976:9) notes that several nineteenth-century mathematicians regarded Euler’s problem about the relation between vertices, edges, and faces in convex polyhedra to be settled beyond doubt by the Cauchy-inspired argument.}. Discovering that cubes, tetrahedra and dodecahedra satisfy an algebraically similar formula than a manipulated version of \eqref{e1} is certainly suggestive. Following Hadamard (1954), we may call the algebraic resemblance a ‘hooking’: a hint. However, we cannot but feel some dissatisfaction with the reasoning. Even if manipulating \eqref{e1} into \eqref{e5} was somehow illuminating, it is hard to imagine what underlying explanation behind the algebraic similarity there might be (besides that it is merely a coincidence). The algebraic description does not provide us with a way of \emph{reasoning} about the regularity observed for cubes, tetrahedra, and dodecahedra, that would accompany the mind in the process of finding out why all polyhedra (and not merely the observed ones) obey the alternating sum regularity expressed by \eqref{e4}.  

By contrast, the Cauchy-inspired geometrical reasoning is stronger. Even if we knew no topology, we may sense that something deep about the geometry of convex polyhedra is invoked. In particular, we cannot but feel that the imaginary operations involved in the reasoning could be performed on any smashed polyhedron that one could conceive: in compressing it, some initially unconnected vertices will suddenly be found on the same surface and will thus be potential vertices for (roughly) triangular figures. In each case, it is not clear what could go wrong and falsify equation \eqref{e4}. Moreover, we know from the setup of Euler’s problem that the length, area, etc. of a figure’s elements are irrelevant to the problem. Accordingly, even if we could not check that \eqref{e4} holds for all cases, it seems that the support that the geometrical reasoning provides to \eqref{e4} is high - higher, indeed, than the one that results from the algebraic reasoning. In Hadamard’s (1954) terms, in coming to learn that smashed cubes, tetrahedra, and dodecahedra obey equation \eqref{e4} we seem to have reached a “relay-result”- a step forward in the direction of the proof.

In summary, the example of Euler’s characteristic suggests that there is an intuitive difference in strength among analogical inferences in mathematics. Some, such as Polya’s algebraic reasoning, seem to work on a mostly heuristic level, in that they mostly aim to make some conjectures salient. Others, instead, such as the Cauchy-inspired geometrical reasoning, are stronger, to the point that we naturally form an expectation about the truth of a conjecture on their basis. In the former case, the trained mathematician tends to see the analogy as a ‘hooking’: a mere suggestion that some interesting regularity might hold. In the latter, the analogy is a ‘relay-result’. Of course, such a distinction is not always sharp. However, we think that it reflects an important phenomenological difference. In the next subsection, we will consider a more graded case, where telling apart a mere analogical ‘hooking’ from a ‘relay-result’ is harder.
\subsection*{\normalsize{2.2 The Basel Problem}}
In 1650, mathematician Pietro Mengoli posed the following problem: which is the precise summation of the reciprocals of the squares of natural numbers? The question can be formalized by asking the precise sum of the following series:
\begin{equation*}
    \sum_{n=1}^{\infty} \frac{1}{n^2}=1+\frac{1}{4}+\frac{1}{9}+\cdots.
\end{equation*}
The first attempted solution is due to Euler in 1735. Below we will present an argument for this solution, based on two subsequent analogies (cf. Bartha 2009:157--159). The first analogy concerns the concept of the roots of a function.  Let us consider the roots of the function:
\begin{equation*}
    \frac{\sin(x)}{x}.
\end{equation*}
These roots are $k\cdot \pi$, where $k$ is an integer number different from $0$. It is possible to construct a infinite product having the same roots:
\begin{equation*}
    \left (1+\frac{x}{\pi}\right )\left (1-\frac{x}{\pi}\right )\left (1+\frac{x}{2\pi}\right )\left (1-\frac{x}{\pi}\right )\cdots.
\end{equation*}

From the fact that they possess the same roots, one could suppose that the two expressions describe the same function. (Notice that this reasoning does not hold in general. However, it points to a deep commonality among the expressions, which will find definition in the theory of analytical functions and will lead to the \emph{Weierstrass factorization theorem} a century later). 

Now we highlight the similarity in form between the expression of power series:
\begin{equation*}
    \sum_{n=0}^{\infty}a_nx^n=a_0+a_1x+a_2x^2+\cdots,
\end{equation*}
and that of polynomials:
\begin{equation*}
    \sum_{n=0}^{N}a_nx^n=a_0+a_1x+a_2x^2+\cdots+a_Nx^N,
\end{equation*}
where $N$ is an integer. Moreover, for a polynomial with $n$ roots $(\alpha_1,\dots,\alpha_n)$, one can show that:
\begin{equation*}
a_0+a_1x+\cdots+a_nx^n=a_0\left(1-\frac{x}{\alpha_1}\right)\cdots\left(1-\frac{x}{\alpha_n}\right)    .
\end{equation*}

Considering the function $\sin(x)/x$ and its power series (expressed as a Taylor series around zero):
\begin{equation*}
 \frac{\sin(x)}{x}=1-\frac{x^2}{3!}+\frac{x^4}{5!}+\cdots, 
\end{equation*}

one should obtain, by analogy, the equality:
\begin{equation*}
\left(1+\frac{x}{\pi}\right)\left(1-\frac{x}{\pi}\right)\left(1+\frac{x}{2\pi}\right)\left(1-\frac{x}{2\pi}\right)\cdots= \frac{\sin(x)}{x}=1-\frac{x^2}{3!}+\frac{x^4}{5!}+\cdots  .
\end{equation*}
By means of a second analogy with the polynomial case, one could generalize the following finite dimensional formula to the infinite case:
\begin{equation*}
    a_1=-a_0\left(\frac{1}{\alpha_1}+\cdots +\frac{1}{\alpha_n}\right ),
\end{equation*}
and (from the conclusion of the first analogy) eventually obtain the relation:
\begin{equation*}
    -\frac{1}{3!}=(-1)\left(\frac{1}{\pi^2}+\frac{1}{4\pi^2}+\cdots\right),
\end{equation*}
which finally after trivial algebraic steps, gives the answer to Basel problem: 
\begin{equation*}
    \sum_{n=1}^{\infty}\frac{1}{n^2}=\frac{\pi^2}{6}.
\end{equation*}

\begin{figure}[ht]
\centering
\includegraphics[scale=0.35]{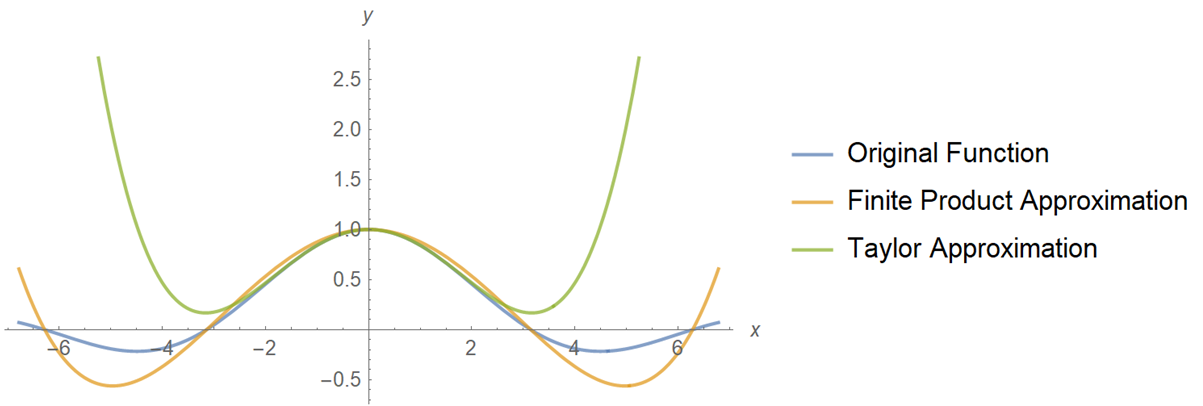}
\caption*{Figure 4. The three curves represent the original function $\sin(x)/x$, the approximation given by the finite product involving the first four roots of the original function, and its fourth order Taylor approximation (actually, the Taylor series of $\sin(x)$ divided by $x$).}
\end{figure}

It turns out that the solution arrived at by the two subsequent analogical inferences is correct.\footnote{Historically, it appears as though Euler boldly took the argument above to constitute a demonstration. Our adaptation brings out the fact that the argument was, in fact, an analogical one - inductive and not deductive.} But how strong was the argument for that solution? The intuitions appear to be less clear. A source of special difficulty is that the inference proceeds to extend an operation (addition) that is known to have meaning in the finite case (polynomials) to a domain of infinities (power series). However, some relevant mathematical properties, such as additivity, may fail to have a clear application in the infinite case. What is at stake, in inferences of this kind, is our entitlement to export an operation. For this reason, the basis for the analogical inference might seem unstable.

Some considerations about the analogical inferences in each case can still be made. The first reasoning, which infers sameness of function from sameness of roots, may feel stronger. Even though the Weierstrass theorem was not known, Euler may have expected that, for at least functions that involve \emph{sine} and \emph{cosine} operations, the inference to the identity of the functions was justified. Considering the matter geometrically, it is plausible that, for at least some \emph{well-behaving} functions, a finite product of linear terms would turn out to be an approximation of the original function - roughly in the same way as Taylor series of finite order approximate the same function (see Fig. 4 for an illustration of the underlying geometrical intuition). These considerations add strength to the mathematical conjecture whereby it is possible to re-describe the original function for polynomials as an infinite product of linear terms depending on its roots. 

The second analogy, which infers the infinite dimensional formula for power series from the finite dimensional formula for polynomials, is instead less convincing. It is harder to avoid the impression that the algebraic structure that is proper polynomials is unlikely to extend so easily to power series. No geometrical intuition of the kind that supports the first analogy appears to be present in this second case. Overall, then, the ‘plausibility’ that the analogy provides to the conjecture seems closer to that of heuristic support, as in ‘hookings’, rather than a ‘relay-result’.

We have now considered two cases of analogical reasoning in mathematics. Our contention has been that there exists an important phenomenological difference in mathematical practice between what we have called analogical ‘hookings’ and ‘relay-results’. While the distinction is not always sharp, with some effort we can often uncover reasons for why some inferences may feel stronger than others to the working mathematician. Can we now identify a general epistemological account that can systematically vindicate and explain the emergence of the differential phenomenology just noted? The next section considers two prominent options.

\section{Theory}
\label{Sec3}
What makes a noted element of analogy a \emph{relay-result} - something that justifies increasing our trust in a mathematical conjecture - as opposed to a mere \emph{hooking}? And what makes an analogical relay-result \emph{more telling} than another? The attempt to provide informative answers to general questions of this kind has traditionally been the province of a philosophical theory of analogical reasoning. The principal object of such a theory are evaluative criteria for analogical arguments. By an ‘analogical argument’ we mean a linguistic object in premise-conclusion form as follows:
\begin{center}
\begin{tabular}{ ccc } 
$A$ is similar to $B$ in respects $P_1,\dots,P_n$ \\
$B$ has some further feature $P_m$ \\[0.25cm]  
\hdashline
\\[-0.25cm]
Therefore, $A$ also has (an analogue of) feature $P_m$.
\end{tabular}
\end{center}

Assuming that the resemblance that a given relay-result or hooking express can be given linguistic form, the original epistemological question about the difference between an analogical relay-result and a mere hooking can then be immediately rephrased as: what makes a given argument from some similarity in the premises capable of providing some \emph{inductive support} to the occurrence of the similarity that stands as the conclusion? This would tell us which kinds of similarities are considered capable of rationally increasing our trust in a conjecture. Relatedly, the question about the difference between a given relay-result and another can be rephrased as: what makes an argument from some similarity in the premises capable of \emph{stronger} support to a given conclusion than another such argument with the same conclusion but different premises?

A prominent answer to these questions, developed specifically for the empirical sciences, is due to Hesse (1963). On her account, the evaluation of an analogical argument in science occurs along two distinct dimensions, namely the ‘horizontal’ and the ‘vertical’ dimensions. The former is concerned with ensuring that the similarities in the argument are genuine and not gerrymandered, i.e., artificially introduced to inflate an otherwise implausible inference. To this effect, Hesse imposes a \emph{materiality} requirement on any inductively strong analogical argument:\footnote{Hesse’s ‘material similarity’ is not synonymous with similarity in physical composition. For instance, Hesse (1963:60) considered material analogies between light and sound, despite the fact that the two are physically very different.}
\begin{quote}
\emph{Materiality}: the similarities mentioned in the analogical argument must be ‘material’ or ‘pre-theoretic’, i.e., they must be cases of sharing of features whose significance in a given context can be recognized before and independently of the analogical argument.
\end{quote}

In addition to the genuineness of the similarities, Hesse argues that the evaluation of an analogical argument involves the assessment of the relevance of the similarities mentioned in the premises to those mentioned in the conclusion. This is the vertical aspect of the evaluation. It is expressed in the following further requirement on inductively strong analogical arguments:
\begin{quote}
\emph{Causality}: the similarities in the premises must be relevant to the occurrence of the similarity mentioned in the conclusion: more specifically, it must be a serious epistemic possibility, which may be further supported by the known similarities, that roughly the same causal connections (stronger than statistical correlations) that hold between the properties of the source also hold between the known and the predicted properties of the target.
\end{quote}

Hesse’s final requirement on inductively strong analogical arguments, which pertains to both the horizontal and the vertical dimension of the analogy, can be stated briefly as follows:
\begin{quote}
\emph{No-Essential-Difference}: There must be no known critical difference, i.e., a difference in the properties and relations of source and target whose occurrence would directly and in the first instance undermine the strength of the analogical argument.
\end{quote}

This final clause serves to ensure that the known differences between the analogy’s source and target are at most secondary ones, and thus incapable of significantly affecting the argument’s conclusion. When \emph{Materiality}, \emph{Causality}, and \emph{No-Essential-Difference} are all satisfied, the result according to Hesse is an analogical argument capable of (some degree of) inductive support to its conclusion. In this way, Hesse’s account promises to vindicate the phenomenologically salient difference between merely heuristic and confirmatory analogical arguments in scientific inquiry.

The recognition of horizontal and vertical dimensions of evaluation also yields a natural answer to the question of when an analogical argument is capable of greater inductive support than another: viz., when the similarities along the two dimensions are more numerous and striking. For instance, analogical inferences about a given human population’s response to a drug are obviously stronger when the model is some other human population on which the drug has been tested, rather than (say) a population of mice which has received the same drug. In practice, Hesse denies that anything more than this intuitive account of degree of analogical support is either possible or necessary: as she writes, “it is difficult to imagine a realistic situation of choice between models where the choice obviously hangs on a judgment about which model is more similar to the [target]. [...] In the first place, comparison of such ‘degrees’ is all but impossible in any practical case... [Moreover,] other considerations will often enter at this stage” (1963:115), such as the simplicity and past success of the model, which is not directly related to its similarity.

An outstanding issue with Hesse’s framework has to do with whether it can be naturally extended to the mathematical domain. Bartha (2009, 2019) has raised doubts over this possibility: “mathematics”, he writes, is “a field in which material analogy, in Hesse’s sense, plays no role at all” (2009:43). In fact, Bartha employs the case of mathematics as part of an argument for a new account of analogical reasoning, the \emph{articulation model}, that promises to extend to both empirical and mathematical sciences. The main point of divergence concerns Hesse’s \emph{Materiality} condition, which Bartha finds exceedingly restrictive of the similarities that can be used legitimately in analogical reasoning. On his more liberal model, the acceptability of an analogical inference mainly depends on the existence of a clear relation between the similarities that figure in the premises of the analogical argument and those that figure in its conclusion - what Hesse would have called the \emph{vertical} dimension of the analogy. Furthermore, Bartha countenances a much broader spectrum of acceptable vertical relations than Hesse, which includes not only causal relations but also statistical correlations and deductive entailments.

With regards to mathematics, Bartha (2009) claims that the relevant vertical relation for an analogical argument is given by the relation of deductive entailment that connects some background mathematical assumptions (e.g., set-theoretic, geometric, arithmetic, etc.) to the property of the source domain that is projected to hold in the target. For instance, the vertical relation in the inference from polynomials to power series in the Basel problem is the proof that:
\begin{equation*}
a_0+a_1x+\cdots+a_nx^n=a_0\left(1-\frac{x}{\alpha_1}\right)\cdots\left(1-\frac{x}{\alpha_n}\right).
\end{equation*}

This view yields the following as a first necessary condition for inductive support by analogy:
\begin{quote}
\emph{Prior Association}: There must be a clear relation of proof, in the source, from the properties $P_1,\dots,P_n$ (i.e., those that are known to be shared with the target) to $P_m$ (i.e., the property that the source possesses but is merely predicted to hold in the target).
\end{quote}

Once the source proof has been articulated, Bartha claims the acceptability of an analogical argument in mathematics depends solely on the following two additional requirements:
\begin{quote}
\emph{Overlap}: Some explicit assumption in the source proof must correspond (admissibly) to a fact known to be true in the target domain.
\end{quote}
\begin{quote}
\emph{No-Critical-Difference}: No explicit assumption in the proof can correspond (admissibly) to something known to be false in the target domain. 
\end{quote}
Under appropriate specifications of the admissibility clauses (to be discussed in the next sections), Bartha’s account yields the verdict that Euler’s analogical argument in the Basel problem was plausible - and similarly for several other cases of analogical inferences in pure mathematics. Furthermore, Bartha defends roughly the following proposal with regards to the relative strength of an analogical argument compared to another (cf. Bartha 2009:173):
\begin{quote}
\emph{Ranking}: Analogical argument $A_2$ is superior to analogical argument $A_1$ for the same conclusion $C$ if and only if: (a) the positive analogy (i.e., the class of known similarities) in $A_2$ includes all critical factors (i.e., all factors relevant to $C$) in the positive analogy of $A_1$, but not vice versa; or (b) $A_2$ and $A_1$ are on a par with regards to condition (a), but in $A_2$ the the source possesses property $P_2$ and in A1 the source possesses property $P_1$, and $P_2$ is closer than $P_1$ to some property $P_3$ that the target possesses (by some natural parameterization). 
\end{quote}
We care to emphasize that neither Hesse’s nor Bartha’s views are intended as answers to the problem of justification. The accounts for, respectively, science and mathematics that each provides purport to merely describe norms and are compatible with a number of approaches to the issue of their justification (including broadly ‘externalist’ ones as defended by Norton 2020). 

Our aim in the following sections is to show that, for all the merits of his account, Bartha was too quick to dismiss Hesse’s framework for material analogy as the basis for a plausible account of analogical reasoning in mathematics. As we will argue, the articulation model often fails to account for the perceived difference in strength found among analogical inferences in mathematics. Specifically, in the next section we will discuss how substituting \emph{Materiality} with \emph{Overlap} results in a clear loss of discerning power. This will prepare the ground for the discussion of sections five and six, where we will argue that the articulation model misses some other important insights of Hesse’s account related to the idea of relevance of similarities and the essentiality of differences. This will once again illustrate how closer adherence to Hesse’s account results in an improved understanding of inductive support by analogy in mathematics. 

\section{Materiality}
\label{Sec4}
A central tenet of Hesse’s account is that no analogical argument can provide genuine inductive support to a conclusion unless there exist similarities in ‘material’ or ‘pre-theoretic’ respects between source and target. This condition immediately rules out as illegitimate analogical inferences based on gerrymandered similarities, such as those of the ‘grue’ variety made famous by Goodman (1954), whenever artificially introduced to inflate an otherwise implausible inference. The same condition also rules out inferences from purely ‘formal’ analogies. One of Hesse’s examples is the analogy between father and child, on the one hand, and state and citizen, on the other. Inferring that citizens owe respect and obedience to the state just as children do towards their fathers is a bad inference, on Hesse’s view, because “there does not seem to be any horizontal relation of similarity between the terms, except in virtue of the fact that the two pairs are related by the same vertical relation” (1963:63). In other words, the similarities between citizen and son do not exist \emph{before} and \emph{independently} of the vertical relation made salient by the analogical argument, such as the relation “provider-for” (1963:63) that links a father to his child. 

Bartha (2009) rejects Hesse’s position on grounds that “many creative analogies in mathematics involve manipulations that are hard to distinguish from specious moves” (111). In other words, the concern is that Hesse’s condition imposes too narrow a limit to what analogical inferences in mathematics can be plausible. His strategy for liberalizing Hesse’s account goes as follows. First, Bartha defines a syntactic notion of similarity for mathematics, distinguishing algebraic and geometric similarity. The former has a straightforward definition: in short, a is algebraically similar to b just in case the respective expressions satisfy some common description . As for geometric similarity, Bartha defends the following criterion: 
\begin{quote}
“Expressions $F$ and $F^*$ are geometrically similar if one of the following three cases obtains: 
\begin{itemize}
    \item[(a)] $F= F(m)$ and $F^*= F(n)$, where $F(k)$ is an expression parameterized by a positive integer $k$ in some range including $m$ and $n$;
    \item[(b)] $F=F(\alpha)$ and $F^*= F(\beta)$, where $F(\tau)$ is an expression parameterized by a real number $\tau$ in some range including $\alpha$ and $\beta$;
    \item[(c)] $F= F(x_1,\dots,x_m)$ and $F^*= F(x_1,\dots,x_n)$, where $F(x_1,\dots,x_k)$ is an expression with $k$ arguments, defined for a range of values of $k$ including $m$ and $n$.” (Bartha 2009:166)
\end{itemize}
\end{quote}
Second, to further restrict the range of similarities that may be used in analogical inference, Bartha defines the notion of an \emph{admissible} similarity. Whereas “all algebraic similarities count as admissible” (2009:165) for Bartha, not all geometric similarities are admissible. Indeed:
\begin{quote}
“For a geometric similarity to count as admissible, the relevant relations and functions should be expressed using standard representation, unless some justification \emph{internal} to the domain can be given for a nonstandard representation. More specifically, any novel representation should be justified in terms of the proof that is the [vertical relation] for the analogical argument.” (Bartha 2009:169).
\end{quote}
The upshot is that ‘admissible similarity’ is proof-dependent. Clearly, this stands in stark contrast with Hesse’s account of analogical inference, according to which inductively significant similarities must be \emph{independent} of the vertical relations introduced by the analogical argument.

As appealing as Bartha’s arguments might seem, the examples below show that abandoning \emph{Materiality} was likely the wrong move: the passage to ‘admissible similarity’ results in a loss of discriminatory power for an account of analogical reasoning in mathematics. Let’s take a look.
\subsection*{\normalsize{Example 4.1. Too many cooks}}
For a start, we note that the account licenses ‘admissible similarities’ that, intuitively, should not be regarded as plausible basis for analogical inference. Consider, for instance, the pair $x^1$ and $x^{-1}$. Geometrically speaking, there is no similarity between the two. Yet, the pair satisfies Bartha’s definition of \emph{geometric similarity} (case (b)). Moreover, since the respective expressions employ the standard representation, their ‘similarity’ is also admissible. The same implausibility affects, a fortiori, the definition of admissible \emph{algebraic similarity}, which is even broader: since, for any two sufficiently long algebraic expressions (a) and (b), we can cook up an almost arbitrarily large number of descriptions $\psi_1,\dots,\psi_n$ that (a) and (b) may or may not satisfy, according to Bartha’s definition it seems that any two expressions will be algebraically similar in an absurdly long number of respects and algebraically dissimilar in an equally long number of other respects.

\subsection*{\normalsize{Example 4.2. Area, Area$^*$, and Area$^{**}$}}
Here is another example of the fact that Bartha’s ‘admissible similarity’ is too broad - in other words, that it licenses far too many analogical inferences with no genuine mathematical ground. From the definition of area and perimeter of a rectangle, one can prove that, of all rectangles of fixed perimeter, the square has maximum area. By analogy, one might suppose that, of all rectangular boxes of fixed perimeter, the cube has maximum volume. After all, there are obvious similarities between perimeter in the two-dimensional and in the three-dimensional case, and between area and volume. For instance, it is well-known that for a rectangle of sides $x$ and $y$ the \emph{Area} can be computed by $A=x\cdot y$ and that, for a rectangular box of sides $x$, $y$ and $z$, the\emph{Volume} of is given by $V=x\cdot y \cdot z$. Given that a square maximizes area, it would be rather surprising if its three-dimensional analogue, the cube, did not similarly maximize volume. 

As Bartha (2009) notes, there exist alternative analogical routes to the conclusion that the cube maximizes volume that are significantly weaker than the argument just provided. For instance, one might reason to the same conclusion by first rewriting \emph{Area} as follows:
\begin{equation*}
A^*(x,y)=3x^{2-2}+x^{2-1}y^{2-1}-3\sin^{2-2}(y),
\end{equation*}
This (equivalent) reformulation induce the following expression for the $3$-dimensional case:
\begin{equation*}
V^*(x,y,z)=3x^{3-2}+x^{3-1}y^{3-1}z3-1-3\sin^{3-2}(z).
\end{equation*}
In effect, there is a resemblance between A* and V* as both satisfy the expression:
\begin{equation*}
3x_1^{n-2}+x_1^{n-1}\cdots x_n^{n-1}-3sin^{n-2}(x_n).
\end{equation*}

Knowing that, of all rectangles of a fixed perimeter, the square maximizes $A^*$, one might infer from the resemblance above  that, of all boxes with a fixed perimeter, the cube maximizes $V^*$.

The challenge is to explain the intuitive difference in strength between the two analogical inferences. For a defender of \emph{Materiality}, the challenge is easily met: the resemblance upon which the second argument is based has no mathematical significance before and independently of the argument; no similarly manufactured resemblance appears, conversely, in the first argument. As we have seen, however, Bartha is not persuaded by this explanation: “if we insist upon a ‘standard way’ of doing mathematics,” he writes, “we might never find the innovative reformulation” (2009:156). Accordingly, he insists that the similarity that the second argument introduces is not suitable for inductive inference, not \emph{qua} manufactured, but \emph{qua} inadmissible relative to the source proof. Specifically, he points out that the proof for the conclusion that the square has maximum area exploits a “symmetry between $x$ and $y$” (2009:170). Since this symmetry is lost in the passage to $A^*$ and $V^*$, the second argument fails to satisfy the condition that “any novel representation should be justified in terms of the [source] proof” (2009:169).

The problem is that it is not difficult to come up with a manipulation satisfying the criterion of symmetry but which is again a totally implausible basis for analogical inference. Consider, for instance, the following function for the area of a rectangle of edges $x$ and $y$, called $A^{**}$:
\begin{equation*}
A^{**}(x,y)=x^{2-1}y^{2-1}-\sin^{2-2}(xy).
\end{equation*}
It is trivial to check that such a formula preserves the symmetry between $x$ and $y$ and that it leads to the following generalization for the volume, namely $V^{**}$:
\begin{equation*}
V^{**}(x, y, z)=x^{3-1}y^{3-1}z^{3-1}-\sin^{3-2}(xyz).
\end{equation*}
where $A^{**}$ and $V^{**}$ are both instances of the general (symmetric) formula:
\begin{equation*}
x_1^{n-1}\cdots x_n^{n-1}-  \sin(x_1\cdots x_n).    
\end{equation*}

In summary, Bartha’s idea to count as admissible respects for similarity manufactured representations justified by the source proof is too broad and leads to invidious discriminations. For instance, Bartha regards the analogical inference based on the resemblance between $A^*$ and $V^*$ as implausible (since it does not display the symmetry that is present in the source proof); but the similarly gerrymandered inference drawing upon the resemblance between $A^{**}$ and $V^{**}$ preserves the symmetry exploited in the source proof and should therefore count as plausible on his account. However, it is clear that the two inferences are equally weak. The problem, in our view, stems from the attempt to make the notion of ‘admissible similarity’ dependent upon the vertical relations in an analogical argument. \emph{Materiality} has a clear advantage in this case, since it treats both the resemblance between $A^*$ and $V^*$ and that between $A^{**}$ and $V^{**}$ as equally gerrymandered and therefore equally incapable of underwriting inductive support by analogy.

\subsection*{\normalsize{Example 4.3. Euler characteristic
}}

Is \emph{Materiality} an expression of chauvinism towards the ‘standard representation’? We think not. Our response is in two parts. Through the following case-study, we aim to illustrate the first part of this response. It is that, in many cases, one must keep in mind the distinction between a heuristic and an inductive sense of analogical support for a mathematical conjecture. As defenders of \emph{Materiality}, we may happily concede that an analogical inference fueled by a manufactured resemblance is often a useful tool for mathematical discovery - a ‘hooking’. Accordingly, we agree with Bartha that without manipulation “we may never find the innovative reformulation” (2009:156). However, not all analogical arguments that are useful heuristically are also inductively strong. In fact, one of \emph{Materiality}’s main advantages compared to \emph{Overlap} is precisely its superior discriminating power with regards to inductive support from analogy.

The example of the Euler characteristic discussed earlier (Section 2.1) offers a beautiful illustration of this. As discussed previously, Polya identifies an ingenious route to the Euler characteristics, based on the algebraic resemblance between the known formulae for two- and three-dimensional figures. Heuristically, noticing that both formulae are alternating sums of the respective geometric elements can be suggestive. However, we refrained from granting it the status of an analogical ‘relay-result’. \emph{Materiality} helps articulate the reason behind the perceived difference: in the geometrical context that governs the example, a similarity with respect to the algebraic description for the two- and three-dimensional case arguably has little significance before the argument is introduced. After all, it is not difficult to imagine alternative three-dimensional conjectures whose algebraic resemblance with the two-dimensional case is, from the perspective of the geometrical context of the example, just as arbitrary as the alternating sum resemblance. Accordingly, the latter is not plausibly regarded as sufficient basis for inductive support. 

The alternative, Cauchy-inspired reasoning is a much clearer candidate for satisfying \emph{Materiality} (see fn. \ref{fnLak}). In this case, the similarities that fuel the inferences are with respect to the geometrical features of polyhedra: for instance, the fact that each polyhedron is a single three-dimensional block made up of a number of faces, connected to each other by edges and vertices. These respects of similarity possess significance before and independently of the Cauchy-inspired argument. They are accordingly capable of sustaining the conjecture that any convex polyhedron will behave in the same way as cubes, tetrahedra and dodecahedra under the ‘compressing’ operation that was described earlier; and that, therefore, each of them will satisfy the very same relation between their vertices, edges, and faces that we referred to as the Euler characteristic. Satisfaction of \emph{Materiality} thus helps explain why the Cauchy-inspired argument appears to possess superior inductive strength within the geometrical context of the example.

\subsection*{\normalsize{Example 4.4. Projection-valued measure as a bridge to infinity}}

Another comment about \emph{Materiality} is necessary, which completes our response to the possible accusation of chauvinism towards the ‘standard representation’. Although \emph{Materiality} expresses a form of caution towards analogical inferences based on newly introduced similarities, it does not deem as heuristic \emph{all} instances of such inferences. To illustrate, let’s consider an actual historical case of analogical inference from finite to infinite. (We anticipate here part of a discussion about asymptotic analogies that will continue in section six). In the branch of mathematics known as spectral theory, ‘self-adjoint’ operators play a fundamental role. They are defined as operators in a finite vector space which possess a special property: they are their own ‘adjoints’. Roughly, this is the equivalent in a suitable functional space of the property whereby, by inverting the rows and columns of a given matrix, one obtains exactly the same matrix again. 

The generalization of such operators from finite vector spaces to infinite-dimensional Hilbert spaces is a non-trivial task. Indeed, by the so-called \emph{spectral theorem}, the condition self-adjointness is a necessary and sufficient condition for the diagonalization of the operator in finite vector spaces. By this, we mean that the eigenvalues of the finite matrix, which are all real, can all be made to appear along a diagonal. If we consider self-adjoint operators on an infinite dimensional Hilbert space, one may be led to infer, by analogy, that a similar mathematical result about diagonalization is plausible. However, the direct generalization fails as many infinite operators have zero as their only eigenvalues, which would induce an inconsistent representation by a zero infinite matrix. In other words, we cannot write a self-adjoint operator on an infinite Hilbert space simply via a linear combination of orthogonal projections, as in the finite case.

Historically, the key to the correct generalization was found by means of a partial redescription of the analogy’s main terms. In particular, a sophisticated mathematical tool called \emph{projection-valued} measures was used. In the finite case, we associate an operator to its orthogonal projections. We then consider rewriting self-adjoint operators in the infinite case as the integral (over a subset of the real line called the ‘spectrum’) of the coordinate function (i.e., the infinite dimensional analogue of the finite dimensional projection coordinates). Hence, the integral of the coordinate function varying over the elements of the spectrum of projection-valued measures becomes the infinitary analogue of a sum of the eigenvalues of the orthogonal projections in the finite case. We thus expect that, as in the finite case the spectral theorem guarantees the existence of a decomposition based on a linear combination of the projection operators and the eigenvalues, so also a generalized spectral theorem for infinite dimensions will yield a similar operation on projection-valued measures and coordinate function.

The example aims to illustrate the following lesson. It is sometimes the case that an analogy between one mathematical domain and another requires a re-description of the original terms in which source and target are described. Such recastings are a delicate matter, since one runs the risk of introducing gerrymandered similarities. However, when the redescription is sufficiently motivated by external considerations, as it arguably happens in this example, \emph{Materiality} does not rule out the possibility of strong analogical inferences.\footnote{Even though we like to keep the issues clearly distinct, we note that a similar claim is plausible with regards to the question whether proofs introducing ‘extraneous’ mathematical concepts can be explanatory; cf. Lange (2019:4).} In this case, the operator-valued measures had mathematical significance in the context of measure theory as extensions of the classical real-valued measures. Accordingly, when the analogical inference was put forth, the redescription of the analogy’s main terms as a similarity between projection operators and projection-valued measure could be regarded as legitimate. Therefore, we have a plausible example in which \emph{Materiality} licenses strong analogical inferences in spite of ‘manipulation’.
\medskip

To summarize, in this section we have argued that inductive support by analogy in mathematics obeys \emph{Materiality} rather than Bartha’s more liberal \emph{Overlap}. In the next section, we will defend another central tenet of Hesse’s account as it applies to the domain of mathematics.
\section{Relevance}
\label{Sec5}

One of the central ideas in Hesse’s (1963) account of analogical reasoning is that, in order for the similarities to provide any additional ground to a conclusion, the properties that the target is known to share with the source must be connected non-accidentally to the properties that are merely predicted to hold for the target. Absent such a connection, the analogical inference is weak. From knowing that John is a Hegel scholar, for instance, one does not plausibly infer that his identical twin Jack is a Hegel scholar, too, despite the impressive bodily resemblances between the brothers. Conversely, sometimes even little in the way of resemblance, as when inferring facts about human beings from observations about mice, can be the basis for a strong analogical inference. Hesse’s way of articulating the difference between relevant and irrelevant similarities goes by means of a separate requirement on the vertical relations in the analogy: they must be \emph{causal} relations “in some scientifically acceptable sense” (1963:87), where, in addition, it must be at least a serious epistemic possibility that causal relations “of the same kind” (1963:87) as the source’s also hold between the known and the predicted properties of the target.

As we have seen, Bartha (2009) is keen on insisting on the vertical dimension of the evaluation of analogical arguments. His approach to translating Hesse’s \emph{Causality} (designed specifically for the empirical sciences) to the domain of pure mathematics goes as follows. As a first step, Bartha imposes his \emph{Prior Association}, whereby the properties of the source that are mentioned in the analogical argument must be related to each other by means of a proof. Further, to ensure that the source proof is a suitable candidate for generalization, Bartha puts forward a \emph{No-Critical-Difference} condition, whereby no explicit assumption in the source proof can correspond to a fact known to be false in the target. The combination of Prior Association and the revised \emph{No-Critical-Difference} condition yields Bartha’s account of analogical relevance.

In what follows, our aim is to show that Bartha’s conditions are once again off-target: in a variety of case-studies from mathematical practice, the conditions prove too weak and unselective. As replacement for Bartha’s proposal, we defend the following alternative:

\begin{quote}
\emph{Relevance}: the vertical relation in the source must be some robust mathematical connection; and it must be a serious possibility that the same mathematical connection that holds in the source also obtains between the known and the merely predicted properties of the target.
\end{quote}

\emph{Relevance} is intended to mirror Hesse’s \emph{Causality} as closely as possible. In particular, we appeal to the notion of a ‘robust mathematical connection’ as the mathematical equivalent of ‘causal relations’ in the empirical sciences. An example is the connection between a polyhedron and its ‘smashed’ two-dimensional counterpart: the latter is connected to the former by \emph{being the result of a given operation} (i.e., smashing) on the former. This connection is clearly of a mathematical nature. Indeed, although such a connection is not itself a proof, it is typically \emph{underwritten} by a proof. Moreover, such a connection is ‘robust’ in that it leads reliably to the predicted outcome. In our example, the smashed polyhedron is the inevitable result of the operation of smashing.

As in Hesse’s original account, a prominent role in the \emph{Relevance} condition is played by the requirement that mathematical connections \emph{of the same kind} as the source’s obtained in the target. Although the general formulation is vague, our case-studies will illustrate how this requirement assumes fairly precise meaning once it is evaluated in a specific context of mathematical inquiry.
One final comment is worth adding here. Although we are (alongside Hesse) suspicious that any precisely statable ranking proposal will capture what makes one analogical inference stronger than another, we find that there is at least a defeasible correlation between what we might call the ‘degree of satisfaction’ of \emph{Relevance} and the strength of an analogical inference. More precisely, it is plausible to think that, other things being equal, the greater our confidence that the same mathematical connections that hold in the analogy’s source also hold in the unfamiliar target (e.g., the smashing operation when extended from tetrahedra to polyhedra of greater complexity), the stronger the analogical inference (to, e.g., the Euler characteristic) will be. This suggestion is not far in spirit from the one in Polya’s (1954b) classic discussion on analogy in mathematics, where it is argued that our confidence in the conclusion of an analogical argument is proportional to “the hope that a common ground exists” (27) between a source and a target.\footnote{In the limiting case in which our ‘hope’ becomes certainty of a common ground, the analogical inference is no longer playing an inductive role, since the conclusion follows directly from our knowledge about the target system.} In practice, a commonality in the respective mathematical connections is often defeasible indication that an underlying mathematical theory exists which explains those commonalities.\footnote{Cf. Mac Lane (1986) on ‘abstraction by analogy’, whereby “a visible and strong parallel between two theories raises the suggestion that there should be one underlying… theory sufficient to give the common results” (436).}

We will return to the issue of why the degree of satisfaction of \emph{Relevance} is only defeasible evidence for the superior strength of an analogical argument in section six. Meanwhile, let’s consider some accessible case-studies that illustrate the distinct need for the \emph{Relevance} condition. 

\subsection*{\normalsize{Example 5.1. Dedekind’s ideals}}

Let’s start with an example in which the respective verdicts of our account (based on \emph{Relevance}) and Bartha’s coincide, but where the explanation for the verdicts is different in each case. In 1801, Gauss formally proved the fundamental theorem of arithmetic, stating that:
\begin{quote}
\emph{Every integer greater than $1$ can be represented uniquely as a product of prime numbers, up to the order of the factors.}
\end{quote}
The proof of this theorem involves the so-called Euclid lemma. When a few years later Kummer introduced the notion of ideal numbers, he asked if any analogous theorem held for them. In general, the answer was negative. For instance, the number $14$ admits of a unique factorization in the set $\mathbb{Z}$ of integers, namely $2\cdot 7$. However, if we consider the following extension of $\mathbb{Z}$: 
\begin{equation*}
\mathbb{Z}\left[-\sqrt{5}\right]=\left\{a+b\sqrt{-5} : a,b\in \mathbb{Z}\right\}
\end{equation*} 
we immediately obtain another factorization of $14$ as follows:
\[
\left(3+\sqrt{-5}\right)\cdot \left(3-\sqrt{-5}\right).
\]

In the 1870s, Dedekind proposed to consider an abstract mathematical object, which (inspired by Kummer’s work) he called the \emph{ideal}. According to the Euclid’s lemma:
\begin{quote}
\emph{If a prime p divides the product ab of two integers $a$ and $b$, then $p$ must divide at least one of those integers $a$ and $b$.}
\end{quote}
In Dedekind’s new theory of ideals, the following analogue of Euclid’s lemma holds:
\begin{quote}
\emph{The ideal $P$ in a ring $R$ is prime if and only, if $P\neq R$ and if $A\subseteq P$, then either $A\subseteq P$ or $B\subseteq P$.}
\end{quote}
He then conjectured that, by analogy with the fundamental theorem of arithmetic, the property of unique factorization by primes that holds for integers also holds for ideals into prime ideals. This, it turns out, was a correct conjecture (at least for an important class of rings in number theory, called the ‘Dedekind’ rings). The argument in this case draws from an algebraic analogy.

From Bartha’s viewpoint, Dedekind’s reasoning was strong: since no explicit assumption in the source proof (which involves the Euclid Lemma) plausibly corresponds to something known to be false in the target, by Bartha’s \emph{Prior Association} and \emph{No-Critical-Difference} Dedekind’s argument to the unique factorization of ideals counts as plausible.

Our framework gives the same verdict but for a different reason: the strength of Dedekind’s argument depends on our recognition that it is a serious possibility that the same robust mathematical connection that holds between the properties of integers and that of unique factorization by primes also holds between Dedekind’s ideals and the unique factorization by prime ideals. Such a mathematical connection is underwritten by the proof that every integer greater than 1 is uniquely factorizable into primes, but it is not itself the proof. Indeed, the connection would be in place even if we considered a different proof (employing different assumptions) for the unique prime factorization of integers.  

In summary, we have an example in which the verdicts of Bartha’s account and of our alternative’s coincide, but where the explanation in each case is different. While the convergence of the two accounts in this example might suggest that the difference is merely terminological, the next three case-studies will show that there is a distinct need for the \emph{Relevance} requirement: Bartha’s \emph{Prior Association} and \emph{No-Critical-Difference} conditions are insufficient to tell apart analogical arguments in mathematics that are capable of inductive support (analogical ‘relay-results’) from weaker arguments possessing merely heuristic potential (‘hookings’).

\subsection*{\normalsize{Example 5.2. Punctured ball}}

Here is an example of a geometrical analogy. We define a \emph{punctured set} as a set in which a point is removed. For instance, the interval $(-1,1)$ without the point $0$ is punctured. Furthermore, we say that a set is \emph{disconnected} if it can be written as the union of two disjoint nonempty sets. The \emph{unit ball} in $\mathbb{R}^n$ is the set of all points for which the distance from the origin is less than one. It is fairly trivial to check that, in one dimension, the punctured unit ball identified by the union:
\[
(-1,0) \cup (0,1)
\]
is disconnected in two components. A possible question is the following: is it true that the punctured unit ball in three dimensions is similarly disconnected? The answer, it turns out, is negative: a punctured ball in $\mathbb{R}^3$  (the sphere) is not a disconnected set; only its analogue in $\mathbb{R}$ is.

As Bartha (2009:156) acknowledges, this example poses a challenge to his account. The analogical inference from the one-dimensional to the three-dimensional case passes his criteria of \emph{Prior Association}, \emph{Overlap}, and \emph{No-Critical-Difference}: there is a clear relation of proof, in the source, going from the geometrical definitions to the conclusion that the punctured unit ball in one dimension is disconnected; a unit ball in one dimension and its three-dimensional analogue also satisfy many common descriptions; finally, no explicit assumption in the source proof corresponds (admissibly) to something false in the three-dimensional case. However, Bartha concedes that, upon reflection, the analogical inference “is intuitively implausible” (2009:156). 

Bartha’s response to the challenge invokes his \emph{Ranking proposal} (presented in section three). For a start, one notes that the conjecture fails to hold for the two-dimensional case: the punctured unit ball in $\mathbb{R}^2$ (the circle) is not disconnected. Therefore, we now have a second analogical inference, from the two-dimensional to the three-dimensional case, suggesting that the punctured unit ball in three dimensions is not disconnected. Moreover, by Bartha’s lights, a clear ranking exists between the analogies: “it is natural to think that $\mathbb{R}^2$ would be a better model for $\mathbb{R}^3$ than $\mathbb{R}$ would be” (Bartha 2009:172). Since  “suboptimal arguments may be ignored” (2009:173) when multiple analogies indicate different conclusions about a target, we can reject the analogical inference from the one-dimensional case based on “the existence of a competitor that is manifestly superior” (173). This yields the desired conclusion that the inference was implausible.

However, there are several problems with Bartha’s response. First, it yields the desired conclusion only when one additionally notes that the unit circle is not disconnected. But one might have expected that the account would yield the right verdict about the implausibility of the original analogical inference quite independently of the evidence regarding the two-dimensional case. Secondly, the account yields the right verdict only under the assumption that $\mathbb{R}^2$ is a better model for $\mathbb{R}^3$ than $\mathbb{R}$ is. It is not clear why this should be so: as Bartha acknowledges elsewhere, “sometimes, in thinking about a problem in three dimensions, the one-dimensional case might be a better guide than the two-dimensional case” (2009:156), since the predicted property may hold in odd dimensions and not in even ones. Hence, it seems that Bartha’s account manages to reach the intuitively right verdict only if one inserts \emph{by hand} that $\mathbb{R}^2$ is a better model for $\mathbb{R}^3$ than $\mathbb{R}$ is.

A different approach to the problem relies on the \emph{Relevance} requirement that we propose. Starting with the one-dimensional source, one notes that the unit ball in $\mathbb{R}$ (the interval) is disconnected by a geometrical element of dimension \emph{zero}, namely the point. One way of criticizing the analogical inference to the three-dimensional case, then, is to note that it purports to show that a sphere can be disconnected by a geometrical object of a much lower dimension. The operation of puncturing is therefore unlikely to disconnect a three-dimensional object as it does for a one-dimensional object. Phrased in \emph{Relevance}’s terms, it is implausible to suppose that the same kind of robust mathematical connection that exists among the properties of the source (the connection between being punctured and being disconnected as it applies to the unit ball in $\mathbb{R}$) may also obtain between the known and the merely properties of the three-dimensional target. 

By adopting a graphic representation, one can see immediately that the \emph{dimensional distance} between the geometrical objects involved in the analogical inference is not preserved:

\begin{center}
\begin{tabular}{ |c|c|c|c|c|c| }
\hline
\textbf{Object} & \textbf{Dimension} & \textbf{Object} & \textbf{Dimension} & \textbf{Object} & \textbf{Dimension}
 \\ 
\hline
Interval & 1 & Circle & 2 & Sphere & 3\\
\hline
Point & 0 & Point & 0 & Point & 0 \\
\hline
Dim. Dist. & 1 &  Dim. Dist. & 2 & Dim. Dist. & 3 \\
\hline
\end{tabular}
\end{center}
However, one can fix the issue with the original analogical inference through the following table:
\begin{center}
\begin{tabular}{ |c|c|c|c|c|c| }
\hline
\textbf{Object} & \textbf{Dimension} & \textbf{Object} & \textbf{Dimension} & \textbf{Object} & \textbf{Dimension}
 \\ 
\hline
Interval & 1 & Circle & 2 & Sphere & 3\\
\hline
Point & 0 & \emph{Line} & 1 & \emph{Plane} & 2 \\
\hline
Dim. Dist. & 1 &  Dim. Dist. & 1 & Dim. Dist. & 1 \\
\hline
\end{tabular}
\end{center}

In non-graphic terms, a more plausible conjecture that draws from the analogy with the puncturing of the unit ball in $\mathbb{R}$ has the form: \emph{the unit ball in $\mathbb{R}^3$ is disconnected into two components by a geometrical object of dimension $2$ (also called a hyperplane)}. Indeed, in this case, it is a serious possibility that a mathematical connection of the same kind as the source’s also obtains between the known and the merely predicted properties of the target. We conclude that \emph{Relevance} helps explain the difference between stronger and weaker analogical arguments; by contrast, Bartha’s account reaches the right verdict only by means of \emph{ad hoc} assumptions.

\subsection*{\normalsize{Example 5.3. Centroids from triangle to tetrahedron}}

Let’s consider another example of geometrical analogy among those that Bartha analyzes. The source in this analogy is the geometry of triangles. By the ‘median’ of a triangle, we mean the segment that unites the midpoint of a side to the opposite vertex. The following theorem holds:
\begin{quote}
\emph{The three medians of a triangle intersect in a common point (a ‘barycenter’).}
\end{quote}
The question is whether this fact about triangles supports, by analogy, the conjecture:
\begin{quote}
\emph{The four medians of a tetrahedron similarly intersect in a common point},
\end{quote}
where ‘median’ here means the segment that unites a face’s centroid with the opposite vertex.\footnote{As Polya (1954) points out, other conjectures in solid geometry are made plausible by the analogy with the theorem about triangles above. This is not in itself a problem for a theory such as ours since we do not think that there is such a thing as ‘the’ analogue target for any given source; but see also section six on multiple analogies.}

According to Bartha (2009:110), the strength of the analogical inference from the two-dimensional to the three-dimensional case depends on the source proof that one picks as the ‘prior association’. In one version, the source proof involves Ceva’s theorem and is purely geometric. From the definition of median of a triangle $ABC$, we have:
\[
\frac{AX}{XB}\cdot\frac{BY}{YC}\cdot\frac{CZ}{ZA}=1
\]
where $X$,$Y$, and $Z$ are the three midpoints of the segments $AB$,$BC$, and $CA$ respectively. It is a consequence of Ceva’s theorem that segments $AY$, $BZ$, and $CX$ are concurrent in a point. An alternative source proof employs analytic geometry and the Euclidean scaffold (i.e., a Cartesian coordinate system). The definition of a midpoint in a triangle is:
\[
X=\frac{A+B}2;	\qquad Y=\frac{A+C}2; \qquad Z=\frac{B+C}2.
\]
By considering each median in its analytical form (e.g., $Xt+(1-t)C$ for median $XC$), we calculate that there is a point that lies on each median, corresponding to the value $t=2/3$. 

By Bartha’s lights, the case-study demonstrates an interesting form of proof-sensitivity: on his view, the analogical inference to the three-dimensional case is “not plausible” (2009:110) when we take the demonstration via Ceva’s theorem as the source proof; it is only “plausible” (110) when we take the demonstration via analytic geometry as the source proof. In the latter case, indeed, we can rely on the observation that the centroids of a tetrahedron are:
\begin{align*}
X&=\frac{A+B+C}3;	&&Y=\frac{A+B+D}3;\\
Z&=\frac{B+C+D}3; 	&&W=\frac{A+C+D}3.
\end{align*}
In this way, we can exhibit a resemblance with the definition of midpoints in the two-dimensional case (see above). By contrast, the analogical inference via pure geometry apparently suffers from the fact that “there is no clear three-dimensional analogue of Ceva’s theorem” (110) and therefore Bartha’s “\emph{No-Critical-Difference} condition fails” (110).

Bartha’s attempt at illustrating his thesis of proof-sensitivity by this geometrical example is unconvincing. One issue to note in passing is that it is not clear how the lack of a “clear three-dimensional analogue of Ceva’s theorem” could be a violation of Bartha’s \emph{No-Critical-Difference}. The latter imposes that “no explicit assumption in the source proof can correspond to something known to be false in the target” (2009:159, our emphasis). However, the fact that it is hard to think of a plausible candidate for a three-dimensional analogue of Ceva’s theorem does not mean that there must be none. If we follow the letter of \emph{No-Critical-Difference}, then, one ought to conclude that analogical inference in pure geometry is equally plausible as that in analytic geometry and, hence, that the above is no genuine case of proof-sensitivity.

Quite independently of this first problem, we find that Bartha’s account gets the plausibility judgments about the case-study exactly wrong. The strength of the analogical inference from the two-dimensional to the three-dimensional case depends mainly, on our view, on a purely geometrical intuition: given that the ‘median’ in both the two-dimensional and the three-dimensional case is what unites the ‘center of mass’ (so to speak) of a figure or face to the opposite vertex, it is plausible that the medians in both the triangle and the tetrahedron case will concur in a point. Once we abandon the context of pure geometry and move to representing triangles and tetrahedra in a Cartesian coordinate system, the analogical inference arguably weakens. First, the similarities upon which we rely become thinner. Indeed, we are left merely with a similarity in form between the expressions for midpoints and centroids, namely:
\[
\frac{A+B}2, \	\frac{B+C}2, \ \text{etc.}  \quad \text{resemble} \quad   \frac{A+B+C}3, \	\frac{A+B+D}3, \ \text{etc.}
\]

Second, and perhaps most importantly, abandoning the context of pure geometry results in a loss of generalizability of the vertical relations in the analogy. When considering the problem in purely geometrical terms, it immediately strikes us as a serious epistemic possibility that the passage from triangles to tetrahedra does not affect the geometrical connection between that the triangle’s midpoints (their ‘centers of mass’), medians, and barycenter in a triangle. Indeed, knowing the fact about triangles, it seems reasonable to expect that a similar geometrical connection holds not just for three-dimensional Euclidean spaces, but in spherical and hyperbolic geometry as well. By restricting ourselves to analytic geometry, instead, our trust in the possibility of generalizing connections that we know to hold in that restricted context is eroded. For instance, while in principle it is possible to obtain an analytical generalization of the theorem about triangles in higher dimensions or in alternative spaces, in each case the formalization is far from trivial and, most importantly, it does not arise naturally from the analytic background.\footnote{Among other things, it is necessary to fix a spherical (resp. hyperbolical) coordinate system to obtain the analytical result in spherical (resp. hyperbolical) geometry. This is a non-trivial mathematical exercise. See also the debate (on which we take no stance here) between ‘fusionists’ and ‘purists’ reviewed in detail by Arana and Mancosu (2012).}

In summary, we hope that the centroid example illustrates in an accessible manner why the idea that the strength of an analogical argument depends on the choice of source proof is flawed. On our view, it is a mistake to recommend (as Bartha’s account does) that the analogical inference from triangles to tetrahedra be evaluated within the context of analytic rather than pure geometry, simply because of some property of the \emph{proof} of the theorem in the source. In effect, once we shift to describing the problem in analytic terms, the similarities that the inference relies upon become thinner and the vertical relations are arguably less suitable candidates for projection. On our view, the analogical inference is much stronger when we consider it from the  context of pure geometry that is its natural home. In that case, we are able to more clearly focus on the geometrical reason for why two seemingly different elements, namely median of a triangle and median of a tetrahedron, are both appropriately named ‘medians’; from this underlying intuition, we can plausibly conjecture that in both cases the medians concur to a point. Once again, this is evidence in favor of our account of relevance and against Bartha’s alternative.

\subsection*{\normalsize{Example 5.4. Abelian groups}}
One final example, this time from an algebraic context, will offer a further illustration of the weakness of Bartha’s account of analogical relevance vis-à-vis ours. In group theory, a \emph{group} $G$ is defined as a set closed under an operation (such as multiplication) satisfying three fundamental properties: associativity, existence of an identity element, and existence of the inverses for all the members of $G$. We say that a group is \emph{abelian} if the commutative property holds for any element of the group, i.e., $ab=ba$ for any $a$,$b$ in $G$. A theorem of group theory states the following:
\begin{quote}
\emph{If $G$ is a group such that $(ab)^2=a^2b^2$, for all $a$,$b$ in $G$, then $G$ is abelian.}
\end{quote}
The proof of this theorem requires only the definition of associativity, of exponential operation and the laws of cancellation (left and right) for multiplication. Based on this theorem, one may wonder whether the following conjecture in group theory is plausible by analogy:
\begin{quote}
\emph{If $G^*$ is a group such that $(ab)^3=a^3b^3$, for all $a$,$b$ in $G^*$, then  $G^*$ is abelian.} 
\end{quote}

According to Bartha, the conjecture is made plausible by the algebraic resemblances between the respective definitions of the elements of $G$ and $G^*$, together with the absence of critical differences. In particular, Bartha notes that all the assumptions needed to prove that $G$ is abelian transfer to the case of $G^*$. On his view, this fact is sufficient to ensure that the source proof for the abelian nature of $G$ is appropriately generalizable. At a closer look, however, the matter is not quite as straightforward. By associativity, one can rewrite the respective equalities as follows:
\[
(ab)^2=a^2b^2=aabb \qquad (ab)^3=a^3b^3=aaabbb
\]
While the two expressions might look similar to the outsider, this point of view actually helps highlight the algebraic \emph{differences between} the expressions. For once we have ‘$aaabbb$’, we may easily expect that the assumptions needed to show that $G^*$ is abelian are significantly greater in number than those needed to show that $G$ is abelian. After all, it is reasonable to expect that the left and right cancellation laws employed for $G$ will not suffice to recover the abelian property in the form  $ab=ba$ in $G^*$, but further properties will have to be invoked in addition.  

Here is the lesson that we draw from this example. In some cases, the fact that no explicit assumptions in the source proof corresponds to something \emph{known} to be false in the target domain is not sufficient to show that there is genuine potential for generalization (as demanded by \emph{Relevance}). For instance, although we know that associativity, exponential operation and cancellation laws suffice to prove that the source group $G$ is abelian, we are already in a position to expect that the same properties will not suffice to prove that the target group $G^*$ is abelian. Phrased in \emph{Relevance}’s terms, we are in a position to expect that the same mathematical connection that holds for the source domain does not hold for the target. Hence, the analogical inference based on the resemblance between members of $G$ and $G^*$ is merely heuristic. Once again, this shows that \emph{Prior Association} and \emph{No-Critical Difference} are too weak and unselective. 

\medskip

In summary, in this section we have argued that there is a distinct need for a condition of \emph{Relevance} on inductively strong analogical inferences. Our appeal to the notion of a ‘robust mathematical connection’ as the equivalent of Hesse’s ‘causal relations’ has been clarified by appeal to several case-studies. In the next section, we will defend a final requirement on inductive support by mathematical analogy, again employing various case-studies to illustrate it.

\section{Essential Differences}
\label{Sec6}

The final clause in Hesse’s (1963) account is the condition whereby, for an analogical argument in science to possess inductive strength, there must be no \emph{essential differences} between the analogy’s source and target. To illustrate the need for this extra requirement, Hesse offers the example of early theories of heat (1963:89). Even though modeling heat as a fluid could capture many aspects of its physical behavior, by the mid-nineteenth century it was understood that heat is not a conserved quantity. At that point, to continue using the fluid model for prediction involved a decision as to whether conservation was essential to fluids or merely an accessory. Importantly, this decision can be seen as transcending the issue of whether there were relevant similarities with the target: one may concede that \emph{Materiality} and \emph{Causality} were indeed met. The further question regards the acceptability of adopting an imperfect model for analogical inference. That decision, Hesse argues, rests “partly on judgments of what is causally essential to the model, and partly on the current availability of alternative models” (92). Not surprisingly, the fluid model was abandoned partly in concurrence with the development of thermodynamics.

We believe that a similar clause applies to analogical inference in the domain of mathematics. The need for this extra requirement can be illustrated by those cases in which multiple analogies exist for a given target or when multiple targets can be constructed from a given source. As we have seen, \emph{Materiality} and \emph{Relevance} possess significant discerning power in those cases - in telling apart weaker from stronger analogical inferences. However, this discerning power can be further refined by including an analogue of Hesse’s final clause on inductive support by analogy:
\begin{quote}
\emph{No-Essential-Difference}: There must be no known essential differences between source and target, i.e., a difference as to the source’s essential properties and mathematical connections. 
\end{quote}
As in Hesse’s account, what counts as an ‘essential difference’ for us is a contextual matter: the decision rests partly on overall judgments about what is mathematically essential to the source domain, and partly on the availability (or lack thereof) of alternative models for the same target.

Cases of analogical inferences from finite to infinite mathematical domains offer the clearest illustration of the intended use of our \emph{No-Essential-Difference} condition. For instance, we know from the history of set theory that two features of our pre-theoretic notion of the ‘size’ of a set clash when extended to infinity: the set of natural numbers strictly includes that of even numbers, and yet they are the same size. The birth of modern set theory required a complex decision as to which of the intuitive features of the notion of size to retain - and, accordingly, which analogical extensions of the finitary notion of set to infinity to consider acceptable (cf. Mancosu 2009). According to our \emph{No-Essential-Difference}, decisions of this type do not rest solely on the consideration of relevant similarities and dissimilarities between source and target, but include a further component of evaluation - what one might call an ‘all-things-considered’ judgment of acceptability. That is to say, considering the respects in which the target differs from the source and the conclusions about the target that the analogy suggests, one must come to a decision as to whether the source can be an appropriate model for drawing novel inferences about the target.

Given the complexities involved in this and other historical examples of analogical inferences in mathematics, an in-depth analysis of the notion of essential property and its role in the historical development of mathematical theories will require a separate discussion. (Incidentally, the high context-sensitivity of the subject explains why we are skeptical that any proposals about ranking analogical arguments based on their intrinsic features will be able to account for mathematical practice). In what follows, we will pursue the more modest aim of illustrating the spirit of our \emph{No-Essential-Difference} condition in a series of relatively accessible mathematical case-studies, at the same time showing how its verdicts are more plausible than those that follow from Bartha’s rival \emph{No-Critical-Difference} condition. This discussion will conclude the defense of our novel framework for analogical reasoning in mathematics.

\subsection*{\normalsize{Example 6.1. Euler's formula}}

Let’s start with a case where \emph{No-Essential-Difference} is plausibly satisfied and the analogical inference is inductively strong. The celebrated Euler formula, in its general form, states that:
\begin{equation*}
e^{ix}=\cos(x)+i\sin(x).
\end{equation*}
It is possible to derive such a formula by comparing the series expansion of the exponential and trigonometric functions. Starting from the definition of complex exponential function, we have:
\begin{equation*}
e^z=1+\frac{z}{1!}+\frac{z^2}{2!}+\frac{z^3}{3!}+\frac{z^4}{4!}+\frac{z^5}{5!}+\cdots.
\end{equation*}

where $z$ is a complex variable. Fixing $z=ix$, and given $i^2=-1$, we immediately obtain:
\begin{equation*}
e^{ix}=1+\frac{ix}{1!}-\frac{x^2}{2!}-\frac{ix^3}{3!}+\frac{x^4}{4!}+\frac{ix^5}{5!}+\cdots.   
\end{equation*}
At this point, we consider the power series formulation of sine and cosine of $x$, namely:
\begin{align*}
\cos(x)&=1-\frac{x^2}{2!}+\frac{x^4}{4!}+\cdots,\\
\sin(x)&=x-\frac{x^3}{3!}+\frac{x^5}{5!}+\cdots.
\end{align*}
By analogy with the domain of finite sums, one conjectures that it is possible to collect the addends of the complex exponential function as follows (note that Euler did not \emph{know} this):
\begin{equation*}
e^{ix}=\left(1-\frac{x^2}{2!}+\frac{x^4}{4!}+\cdots\right)+i \left(\frac{x}{1!}-\frac{x^3}{3!}+\frac{x^5}{5!}+\cdots\right).    
\end{equation*}

Substituting the terms on the right-hand side with the power series of sine and cosine, we obtain:
\begin{equation*}
e^{ix}=\cos(x)+i\sin(x).
\end{equation*}

On our view, this analogical argument from finite to infinite sums is strong. We know that the operation of sum on power series satisfies several properties analogous to finite sums - hence the term ‘sum’ used for both operations. For instance, consider the two power series:
\begin{equation*}
\left[a_0+a_1x+a_2x^2+a_3x^3+\cdots\right] \quad \text{and} \quad  \left[b_0+b_1x+b_2x^2+b_3x^3+\cdots \right],
\end{equation*}
Their sum is defined in a way entirely analogous to the case of finite sums, as simply:
\begin{equation*}
(a_0+b_0)+(a_1+b_1)x+(a_2+b_2)x^2+(a_3+b_3)x^3+\cdots  
\end{equation*}
The inference that the addends of a complex exponential function could be collected in the exact order required to obtain Euler’s formula, just as one would expect in the case of finite sums, is plausible partly because, as required by \emph{No-Essential-Difference}, we do not know of any feature of the ‘$+$’ operation in the infinite case that is radically different from the analogue operation in the finite case, such that we could expect the re-grouping operation performed on the original exponential function to fail. (As a matter of fact, the regrouping operation is not justified in the general case, but only holds when each series is convergent, as it happens in this example).

Bartha’s \emph{No-Critical-Difference} yields the same verdict in this case. The source proof for the analogical argument is the demonstration that, in the finite case:
\begin{equation*}
1+x+x^2+x^3+\cdots+x^n=\left(1+x^2+\cdots+x^{n-1}\right)+\left(x+x^3+\cdots+x^n\right).
\end{equation*}
The proof involves both commutativity and associativity of addition. Since none of these explicit assumptions in the source proof is known to fail in the infinite case, \emph{No-Critical-Difference} yields the conclusion that the analogical inference to Euler’s formula is plausible. However, as the next example will show, Bartha’s condition faces immediate trouble in related examples.

\subsection*{\normalsize{Example 6.2. The sum of all numbers}}
Let’s consider the derivation, due to Srinivasa Ramanujan, that the sum of all positive integers:
\begin{equation*}
\sum_{n=1}^{\infty}n=1+2+3+4+\cdots= -\frac{1}{12}.
\end{equation*}
The proof proposed by Ramanujan goes as follows. Let’s call the sum of all numbers (whatever it might be) ‘$S$’. By stipulation:
\begin{equation*}
S=1+2+3+4+\cdots
\end{equation*}
We then multiply everything by a factor $4$:
\begin{equation*}
4S=4+8+12+16+\cdots
\end{equation*}	
Using the \emph{additive identity law} (which we know to hold for finite sums), we obtain:
\begin{align*}
S &=1+2+3+4+5+\cdots\\
4S &= 0+4+0+8+0+\cdots
\end{align*}
Subtracting the latter equation from the former, we obtain the alternating series:
\begin{equation*}
-3S=1-2+3-4+\cdots
\end{equation*}
The right side of the resulting equation is the power series of the function:
\begin{equation*}
\frac{1}{(1+x)^2}, 
\end{equation*}
for $x=1$. Hence, Ramanujan concluded that:
\begin{equation*}
-3S=\frac{1}{(1+1)^2}=\frac{1}{4},
\end{equation*}
and therefore, by trivial calculations:
\begin{equation*}
S=-\frac{1}{12}.
\end{equation*}

The conclusion is, of course, complete nonsense.\footnote{One way of interpreting Ramanujan’s result is as indirectly revealing important properties of divergent series.} The obvious culprit in the proof is the extension of the \emph{additive identity law}, which is a property of finite sums, to infinite sums. We highlight that the extension seems to be licensed by Bartha’s account. It is not clear what assumption in the source proof for finite sums (such as commutativity and definition of zero) fails to hold in the case of sums with infinitely many elements. At least \emph{prima facie}, then, Bartha’s \emph{No-Critical-Difference} condition commits us to thinking that Ramanujan’s argument to the conjecture that the sum of all natural numbers is equal to $-\frac{1}{12}$ is a strong argument, based on a plausible analogical extension of the additive identity law from finite to infinite sums.\footnote{We note that the inference meets the \emph{No-Critical-Difference} test even in Bartha’s revised formulation on p. 177. The notion of a ‘dense subset’ of the target does not have a clear application in this case, since finite sums are not a dense subset of infinite sums. Even if they were, it is not clear which one of Bartha’s conditions fails.}

Our \emph{No-Essential-Difference} condition diverges sharply on the verdict. According to it, the plausibility of the extension of the additive identity law to infinite sums rests on an all-things-considered judgment about what is mathematically distinctive about finite sums. In the case of Ramanujan, we have a simple deductive argument to the effect that, if the additive identity law holds in the infinite case, the sum of all numbers is $-\frac{1}{12}$. We are allowed to take this fact into account when considering the overall acceptability of the analogical extension. When we do so, the correct all-things-considered judgment is to reject the extension. After all, one feature that seems to be ‘essential’ to finite sums is that, when we add one positive integer to another, the result is a larger positive integer. This feature is lost in the infinite case if it is true, as Ramanujan’s proof suggests, that the sum of all positive integers is a negative number. Based on this, we conclude that the extension of the additive identity law violates \emph{No-Essential-Difference}.

In summary, the case of Ramanujan’s argument illustrates one important application of our requirement of \emph{No-Essential-Difference} to mathematics. Bartha’s \emph{No-Critical-Difference}, by contrast, fails quite dramatically in the example. It is worth clarifying what the lesson is in this case. In the previous section, we have argued that Bartha’s \emph{No-Critical-Difference} cannot substitute a distinctive requirement of \emph{Relevance}, since it cannot guarantee the potential for generalization. The problem identified in this section is stronger: even if a \emph{Relevance} condition were added to Bartha’s account, his \emph{No-Critical-Difference} condition would still not be sufficient to capture the distinction between strong and weak analogical inferences in mathematics. This is because his condition is solely concerned with the features of the source proof that fail to hold in the target domain. As the case of Ramanujan’s proof illustrates, what we need is a more stringent requirement, which considers not only the known differences in the source proof but also the overall acceptability of the analogical inference in light of the source’s ‘essential’ properties.   

\subsection*{\normalsize{Example 6.3.  Back to the Basel Problem}}
Having considered two extreme examples of (respectively) strong and weak analogical inferences, let’s conclude this discussion with a more graded case: Euler’s reasoning to the solution of the Basel problem. (The reader is invited to look back at section two for an intuitive presentation of the case-study). According to Bartha, who considers this example at length, Euler’s “beautiful exploitation of similarities between polynomials and power series” (2009:157) offers an example of a “plausible” (159) analogical argument in mathematics. As anticipated in section two, our judgment is more mixed: while Euler’s first analogical argument, which infers similarity (indeed, sameness) of function from a similarity in the roots, seems stronger to us, the second argument, which infers the similar behavior of the coefficients of power series and polynomials from the similarities in their respective expressions, seems weaker. The \emph{No-Essential-Difference} condition may play a role in explaining these differential judgments.

For a start, we highlight that there is an important difference in \emph{boldness} between the respective conclusions. With the first analogical argument, Euler aims to establish that: 
\begin{equation*}
\frac{\sin(x)}{x} \quad  \text{\emph{‘is’}} \quad \left(1+\frac{x}{\pi}\right)\left(1-\frac{x}{\pi}\right)\left(1+\frac{x}{2\pi}\right)\left(1-\frac{x}{2\pi}\right)\cdots.
\end{equation*}

The conclusion in this case is in analytic geometry, where the \emph{‘is’} need not be one of strict identity. It is sufficient for the purpose of the argument that the infinite product $(1+x/\pi)(1-x/\pi)(1+x/2\pi)(1-x/2\pi)\cdots$ be a close geometrical \emph{approximation} of the function $\sin(x)/x$,  in the sense which is illustrated in Fig. 4. Conversely, Euler’s second analogical argument aims to establish a much bolder conclusion, viz., that the same relation between coefficients, whereby $a_0+a_1x+\cdots+a_nx^n$ \emph{‘is’} (algebraically speaking) $a_0(1-x/\alpha_1)\cdots(1-x/\alpha_n)$ in the finite domain of polynomials also holds in the infinite domain of power series:
\begin{equation*}
a_0+a_1x+a_2x^2+\cdots   \quad \text{\emph{‘is’}} \quad    a_0\left(1-\frac{x}{\alpha_1}\right) \left(1-\frac{x}{\alpha_2}\right)\cdots.
\end{equation*}
 
In this case, the required notion of \emph{‘is’} is significantly more demanding: it is strict identity. 

Because the conclusion of the second argument is much bolder than the first, we can expect the differences between source and target to have greater weight. One obvious suspect is that there is no equivalent of the notion of \emph{degree} for infinite power series. To wit, the operations in virtue of which one obtains the algebraic identity between the sum and the product representation of polynomials (the finite case) is based on the property of the degree of a polynomial. This degree must necessarily be finite in order for the algebraic identity to have any meaning in the polynomial context. Consequently, Euler’s extension of the algebraic identity to the infinite case of power series naturally raises concerns. Even though the conclusion that we obtain by means of this extension is not nearly as intuitively implausible as the one that results from Ramanujan’s proof, one might still feel some resistance at the idea that the passage to power series should be so straightforward. (With hindsight, the fact that power series sums behave similarly to polynomial sums in the specific case that Euler considered is more of an exception than the rule.)

In summary, while Euler’s first analogical argument appears to be strong, the matter is much less straightforward in Euler’s second inference. We trace this felt difficulty back to the fact that it is hard to determine whether the lack of an equivalent of the notion of degree for power series is an essential or a secondary difference. Indeed, we countenance both the pessimistic view, according to which the lack of an equivalent of degree for power series represents an ‘essential’ difference, sufficient to make Euler’s argument merely heuristic, and the optimistic view whereby the difference, though significant, is not ‘essential’. Importantly, we think that this example illustrates the \emph{No-Essential-Difference} at work. This is because, by appealing to this condition, we can provide a plausible explication of the conflicting judgments regarding the strength of Euler’s analogical inference. We believe that this is a perfectly fine achievement for a descriptive theory: our conditions should tell apart strong from weak analogical arguments whenever the distinction is relatively clear; in borderline cases, we can often only do so much as acknowledge the blurriness and pinpoint its origin with the help of the proposed conditions.

\medskip

To summarize, in this section we have defended a new \emph{No-Essential-Difference} condition, constructed in close parallel to Hesse’s final clause on analogical inferences in the empirical sciences. The three cases of analogical inferences from finite to infinite domains that we have considered illustrate pretty clearly the need for this final requirement on inductive support by mathematical analogy. The next section will wrap up and indicate directions for future research. 
\section{Conclusion}
\label{Sec7}
In this paper, we have offered a novel account of analogical reasoning in mathematics. The above discussion can be thought of as paralleling Hesse’s (1963) pioneering defense of the historiographic and epistemological centrality of the notion of material analogy in the empirical sciences. On our view, a closely related notion is necessary to capture the inductive practice in mathematical research. We should re-emphasize that the account proposed in this paper is intended solely as an answer to the descriptive problem of what makes an analogical inference in mathematics strong. Such a problem is ‘descriptive’ not in the sense that it does not concern norms - indeed, it is very much concerned with identifying the implicit norms of strong analogical reasoning - but in the sense that it does not address the question of their justification - of what makes it rational for epistemic agents like us to reason in broad accordance with those norms. Our account is compatible with a number of potential solutions to the latter problem. 

As brief summary of our view, we have defended the necessity of restricting the inductively significant similarities in mathematics in accordance with the following condition:
\begin{quote}
\emph{Materiality}: the similarities mentioned in the analogical argument must be ‘material’ or ‘pre-theoretic’, i.e., they must be cases of sharing of features whose significance in a given context can be recognized before and independently of the analogical argument.
\end{quote}
Moreover, contra Bartha’s (2009) case for softening the conditions for analogical relevance in mathematics, in sections five and six we have argued that the original division of labor that Hesse established between \emph{Causality} and \emph{No-Essential-Difference} must be restored for the domain of mathematics. Accordingly, we endorsed the following additional requirements: 
\begin{quote}
\emph{Relevance}: the vertical relation in the source must be some robust mathematical connection; and it must be a serious possibility that a mathematical connection \emph{of the same kind} as the source’s also obtains between the known and the merely predicted properties of the target.
\end{quote}
\begin{quote}
\emph{No-Essential-Difference}: There must be no known essential differences between source and target, i.e., a difference as to the source’s essential properties and mathematical connections. 
\end{quote}

As a proposal about what makes an analogical inference in mathematics strong, the above account may gain new insights from historical and philosophical work on specific areas or aspects of mathematical research. Examples such as Euler’s reasoning in the Basel problem already indicate that, for instance, what counts as an essential difference in geometry is rather different from what counts as an essential difference in algebra. Further research in the history and present state of mathematics is likely to uncover similar kinds of field-specific features that escape the depth of field of our account. Even so, the wider view provided in this paper remains useful, as it contributes to showing that the specific practices emerging from the different areas of mathematics can be subsumed under a more general and unifying evaluative framework. 

A comment about the vagueness of the proposed conditions is worth adding here. The notions of ‘material similarities’, ‘mathematical connections’, and ‘essential differences’ are clearly vague to some extent; contextual information must typically be filled in to specify their referents. On our view, however, this is an advantage of our proposal. Although our account aims to articulate what makes for a strong or weak analogical argument, it does not - and, we think, it should not - try to say, \emph{of any given} analogical inference in mathematics, whether it is strong or weak. After all, disagreements frequently arise in the mathematical community as to whether some proposed mathematical analogy is deep or superficial. The best we can hope from our account in those cases is not to try to settle the issue on behalf of mathematicians, but rather to help articulate their disagreements. In other words, our account will have accomplished its job if it succeeded in \emph{clarifying} what the opposing parties in the dispute are disagreeing about - presumably, the extensions of the terms ‘material similarities’, ‘mathematical connections’, and ‘essential differences’ - leaving to those very parties the job of settling their disagreement.    

In conclusion, let us indicate a final issue which deserves future attention. While in the above discussion we have relied upon an intuitive notion of an analogical argument providing \emph{inductive support} to a mathematical conjecture, we have systematically refrained from giving any further explications of that notion. It is an interesting and open question whether it can be given more precise meaning and, in particular, whether it can be understood within a fully general probabilistic theory of confirmation for mathematical conjectures. Because the difficulties associated with extending a probabilistic epistemology such as Bayesian confirmation theory to the domain of pure mathematics are considerable (in light of the well-known problem of logical omniscience), we believe that a full discussion of this issue is better left to a separate occasion.

\end{document}